\documentclass[11pt]{article}

\usepackage{amsfonts}
\usepackage{amssymb}
\usepackage{amsmath}

\newtheorem{theorem}{Theorem}
\newtheorem{lemma}[theorem]{Lemma}

\newtheorem{proposition}[theorem]{Proposition}
\newtheorem{conjecture}{Conjecture}

\newenvironment{proof}{\noindent {\it Proof~:}\ }{\ \rule{1mm}{2mm}\medskip}
\newcommand{\beeq}{\begin{eqnarray*}}
\newcommand{\eneq}{\end{eqnarray*}}

\newcommand{\n}{\mathbb N}
\newcommand{\z}{\mathbb Z}
\newcommand{\zp}{{\mathbb Z}/p{\mathbb Z}}

\newcommand{\ox}{\overline{X}}

\newcommand{\ost}{\overline{S+T}}
\newcommand{\oxy}{\overline{X+Y}}

\begin{document}
\title{On the critical pair theory in $\zp$}
\author{Yahya Ould Hamidoune\thanks{Universit\'e Pierre et Marie Curie, Paris.
{\tt yha@ccr.jussieu.fr}}
\and
Oriol Serra\thanks{Universitat Polit\`ecnica de Catalunya, Barcelona.
{\tt oserra@mat.upc.es}}
\and
Gilles Z\'emor\thanks{\'Ecole Nationale Sup\'erieure des
T\'el\'ecommunications, Paris.
{\tt zemor@enst.fr} \newline
Submitted to Acta Arithmetica, september 2003. Revised, june 2005.}}
\date{June 30, 2005}

\maketitle

\begin{abstract}
Let $A, B$ be subsets of $\zp$ such that
$$
|A+B|\le |A|+|B|+1.
$$
We prove that, if $|A|\ge 4$, $|B|\ge 5$, $|A+B|\le p-5$ and $p\ge 53$
then $A$ and $B$ are included in arithmetic progressions with the same
difference and of size $|A|+2$ and $|B|+2$ respectively.
This extends the well-known theorem of Vosper and a recent
result of R{\o}dseth and one of the present authors.
\end{abstract}

\section{Introduction }
\subsection{Context and main result}\label{sec:context}
The
Cauchy-Davenport Theorem \cite{CAU,DAV} states that if $A$ and $B$
are subsets of $\zp$ then
$$
|A+B|\geq \min (p,|A|+|B|-1).
$$
Vosper's Theorem \cite{VOS} solves the related critical pair problem
and states that, if $|A|,|B|\ge 2$, then
  $$|A+B|\geq \min (p-1,|A|+|B|)$$
unless both $A$ and $B$ are   arithmetic
progressions with a common difference. The Cauchy-Davenport Theorem
was generalized to general
abelian groups by several authors including Mann \cite{MANA}
 and Kneser \cite{KN3}.
Kemperman proposed in \cite{KEM} a recursive procedure which
generalizes Vosper's Theorem to all abelian groups.

In the case when the group is $\z$, the analogue of Vosper's theorem
was considerably strengthened. A theorem of Freiman
states \cite{freiman} that if a subset $A$ of integers
is such that $|A+A|$ is sufficiently small then it must be contained in a
short arithmetic progression: specifically, if $|A+A|\leq 2|A|+m$
with $-1\leq m\leq |A|-4$, then $A$ must be contained in an arithmetic
progression of length $|A|+m+1$.

In the case of the sum of two sets, Freiman's Theorem was generalized
by Lev and Smeliansky, a version of it states \cite[Lemma 1]{levsmelianski}

\begin{theorem}\label{thm:ls}
Let $A,B\subset [0,b]$, where $b\in \n$ with $0\in A\cap B$ and $b\in B.$
Assume moreover $\gcd(B)=1$. Then
$$
|A+B|\geq \min (b+|A|,|A\cup (A+b)|+|B|-2).
$$
\end{theorem}

Going back to the case of $\zp$, Freiman proved \cite{FR1},
using trigonometric sums, that
for a subset $A\subset \zp$,  $|A+A|\leq
 2.4|A|$ and $|A|\leq p/35$ hold only if $A$ is contained in a short
 arithmetic progression. Bilu, Lev and Ruzsa \cite{BLR} show how to obtain that,
for $|A|$ small enough, $|A+A|\leq 3|A|-4$ holds only if $A$
 is a subset of a short arithmetic progression. It is also conjectured
that the result should hold without the restriction on the size of
$|A|$. We will refer to this as the  $(3k-4)$-Conjecture in $\zp$.

Going beyond Vosper's Theorem without any undue
restriction on the size of $A$ and $B$
has proved challenging. A recent result of R{\o}dseth and one of the
present authors \cite{HR} is a step in that direction.

For a subset $X$ of an abelian group, let $\ell_r(X)$ denote, if it
exists, the cardinality
of the smallest arithmetic progression with difference $r$ containing~$X$.
We have \cite{HR}:

\begin{theorem}
\label{thm:hamrod}
Let $A, B$ be subsets of $\zp$ with $|A|\ge 3$ and $|B|\ge 4$. If
$$
|A+B|\le |A|+|B|\le p-4,
$$
then $\ell_r(A)\le |A|+1$ and $\ell_r(B)\le |B|+1$ for some $r\in \zp$.
\end{theorem}

The above results suggest the following conjecture.

\begin{conjecture}\label{conj}
Let $m$ be a non-negative integer and let
$A$ and $B$ be subsets of $\zp$ such that
$$
|A+B|\leq |A|+|B|+m\leq p-(m+4).
$$
If $|A|\geq m+3$ and $|B|\geq m+4$ then   there is $r\in \zp$ such that
$$
\ell_r(A)\le |A|+m+1\; \mbox{  and  }\; \ell_r(B)\le |B|+m+1.
$$
\end{conjecture}

If true, the conditions of Conjecture \ref{conj} can not be weakened.
For example,  take  $A=\{0,1,\ldots ,m+1\}\cup \{2m+4+j\}$ for any
positive integer $j$; we have $|2A|=2|A|+m$ for large enough $p$ and $A$
is not contained in  an arithmetic progression of the stated length; now
let $B=\zp\setminus (-2A)$; then it can be easily checked that
$|A+B|= p-|-A|=|B|+|A|+m$, which shows that the condition $|A+B|\le
p-(m+4)$ cannot be removed.

Note that the choice $m=|A|-4$ in  Conjecture \ref{conj} gives  the
$(3k-4)$-conjecture in $\zp$ for sets $A$ with $|A|\le (p-1)/4$. The
best known result due to Green and  Ruzsa \cite{gr} implies only the
validity of this $(3k-4)$--conjecture when $|A|<10^{-180}p.$

The case $m=0$ of Conjecture \ref{conj} is Theorem
\ref{thm:hamrod}. In the present paper we extend Theorem
\ref{thm:hamrod} to the case $m=1$.
More precisely we prove:

\begin{theorem}\label{thm:main} Let $A, B$ be subsets of $\zp$ with
$|A|\ge 4$ and $|B|\ge 5$. If $p\geq 53$ and
$$
|A+B|\le |A|+|B|+1\le p-5,
$$
then there is $r\in \zp$ such that
$$
\ell_r(A)\le |A|+2 \mbox{ and } \ell_r(B)\le |B|+2.
$$
\end{theorem}

\subsection{Methodology: isoperimetric tools and outline of the paper}

 Let $G$ be an abelian group and let $B\subset G$ be a generating subset
 of $G$ such that $0\in B$. Let $k$ be a positive integer:
 the $k$-th {\it isoperimetric number} of $B$ is
$$
\kappa _k (B)=\min  \{
|X+B|-|X|\  \Big| \ \
|X|\geq k \ {\rm and }\ |X+B|\leq |G|-k\},
$$
where $\min \emptyset =|G|$, by convention. A subset $X$ achieving
the above minimum is called a $k$-{\it fragment} of $B$.
A $k$-fragment with minimal cardinality is called a $k$-{\it atom}.

In the context of additive problems,
atoms were first introduced in \cite{HMAN}
and have since proved useful tools in critical pair
theory, see e.g. \cite{HALG,HWAR,HACTA}.

It is proved in \cite{HMAN} that for any abelian group $G$ and
any subset $B\subset G$,
a $1$-atom of $B$ containing $0$ is a subgroup. This result implies easily
Mann's generalization of the Cauchy-Davenport Theorem.
The structure of $2$-atoms is  more difficult
to describe. For $\kappa _2 (B)\leq |B|$, it is proved
in \cite[Theorem 6.2]{HACTA} that a  $2$-atom of $B$ containing $0$ with size $\geq 3$
is a  subgroup. A critical pair theory is deduced from this description
in \cite{HACTA}.
$2$-atoms were again used in \cite{HR} as an essential tool in the proof of
Theorem \ref{thm:hamrod}.
The structure of $2$-atoms in $\zp$ was further studied
by the authors of \cite{SZ}. As a consequence
it is shown that every set $B\subset \zp$ is the union of
$h=\kappa _2(B)-|B|+2$ arithmetic
progressions with the same difference provided that $|B|\le p-(h+2)^2/2$.
This could be seen as a partial result on the way to Conjecture~\ref{conj}:
to cover the remaining ground one would need in particular to bound the gaps
between consecutive arithmetic progressions.

In the present paper we depart from previous work by making use
for the first time of $k$-atoms for $k>2$. We envisage the following
isoperimetric method to deal with Conjecture \ref{conj}: it comes
in three steps.

\begin{itemize}
\item  {\bf Step 1 }: Prove a special case of Conjecture \ref{conj}
assuming that one of the two sets $A$ and $B$ is already known
to be contained in a short arithmetic progression.
\item {\bf Step 2 } : Replace  the hypotheses of Conjecture \ref{conj}
by the weaker conditions $\kappa_{m+4}(A)\leq |A|+m$ and
$\kappa_{m+3}(B)\leq |B|+m$.
 \item {\bf Step 3 } : Study an $(m+3)$-atom $K$ of $A$  and an $(m+4)$-atom
$L$ of $K$. 
It is
enough to prove the result for $L$, since a repeated application of
Step 1 allows one to recover the structure of $A$ and $B$.
\end{itemize}

Before trying this approach for large values of $m$, we think that
better understanding of the structure of the $k$-atoms is required.
Indeed, since Conjecture
\ref{conj} contains the longstanding $(3k-4)$--conjecture for $\zp$, one may
suspect that the above program will not be without technical
difficulties. However, to check the soundness of the proposed method
 we shall try it out in this paper by considering the
first open case $m=1$ of Conjecture \ref{conj}.

The paper is organised as follows. In Section \ref{sec:size} we
give some preliminary results which allow us to give a first bound
on the size of $k$-atoms of a set $B$ for $k\le |B|$. In Section
\ref{sec:comp} we show that if a pair of sets satisfy the
conditions of Theorem \ref{thm:main} and one of the two sets is
contained in a short arithmetic progression then so is the second
one. We  next show in Section \ref{sec:sizebis} that $4$-atoms of
a set $B$ with $\kappa_5(B)\le |B|+1$ have cardinality $4$ and
that when $|B|=4$ then the $5$-atoms of $B$ have cardinality $5$.
Section \ref{sec:small} presents the remaining ingredient of the
proof which consists of saying that, when the two sets in Theorem
\ref{thm:main} are small,  they are indeed contained in a short
arithmetic progression. Finally, the proof of Theorem
\ref{thm:main} is completed in Section \ref{sec:main}.

\section{A first bound on the size of atoms}\label{sec:size}

In this section we introduce notation and preliminary results which will be
used  throughout the paper.
We also derive a preliminary bound on the size of $k$-atoms that we will need.

Two subsets $X$ and $Y$ of $\zp$ will be said to be {\it equivalent}
if there is $u\neq 0$ and $v$ such that $Y=u\cdot X+v$, where $u\cdot X=\{ux, \; x\in X\}$.
Note that the $k$-isoperimetric numbers of a set  and the fact of being a
$k$-fragment or a $k$-atom  are invariant properties under translations and
automorphisms of the additive group $\zp$.
Therefore there is no loss of generality in considering equivalent sets.

We recall known results that we shall use.

\begin{lemma}[\cite{HALG}, Lemma 2.4]
Let  $B$  be a  subset of a finite
 abelian group $G$.
Let  $F$  be a
  $k$-fragment of $B$ and $a\in G$. Then $a-F$ and   $G\setminus (F+B)$
   are $k$-fragments of
    $-B.$
Moreover     $\kappa _k(-B)=\kappa _k(B) .$
\label{dual}
 \end{lemma}

The following is a particularly useful property of $k$--atoms.

\begin{theorem}[\cite{HALG}, Proposition 2.5]
 Let  $B$  be a  subset of a finite abelian group $G$.
Let $M$  be a    $k$-atom of $B$. Let  $F$  be a
  $k$-fragment of $B$ such that  $ M \not\subset F $.
  Then $|M \cap F|\leq k-1 .$
\label{intersection}

 \end{theorem}

Throughout the paper we use the following notation.
Let   $X$ and $Y$ be   subsets of  $\zp$. For each integer $i\ge 0$
we introduce the set $N_i(X,Y)$ defined by
\beeq
&&N_0(X,Y) = X, \hspace{9mm}N_1(X,Y)=(X+Y)\setminus X\\
&&N_i(X,Y) =(X+iY)\setminus (X+(i-1)Y),\; i\geq 2,
\eneq
where $iY=\underbrace{Y+\cdots +Y}_{i}$. We simply write $N_i$
when the reference to sets $X$ and $Y$ is clear from the context.

For a subset  $U$   of $Y$ and $i\ge 1$ such that
$N_i\neq\emptyset$, we denote by  $N_i^{U}$ the set of elements
$z\in N_i$ such that $z-U\subset N_{i-1}$ and $(z-(Y\setminus
U))\cap N_{i-1} = \emptyset$. We also write
$$
N^{\subseteq U}_i=\bigcup_{V\subseteq U}N^V_i.
$$

\begin{lemma}\label{lem:NU}
Let $X,Y\subset\zp$, $U\subset Y$ and $i\ge 1$.
With the notation just introduced, if $N_{i+1}^U\neq \emptyset$ then
$$
N_{i+1}^U -U\subset N^{\subseteq U}_{i}.
$$
In particular,
$$
|N_{i+1}^U |\le   |N^{\subseteq U}_{i} |-|U|+1.
$$
\label{eq1}
\end{lemma}

\begin{proof}
Let $z\in N_{i+1}^U $, $u\in U$ and $z'=z-u\in N_{i} $.
Then $z'\in N^V_{i} $ for some subset $V$ of $Y$.
But, for any $v\in V$, we have  $z-v = z'-v+u \in N_j $
for some $j<i+1$. Since $z\in N_{i+1} $ we must have $j=i$~:
this implies $V\subset U$.  In particular, if $N_{i+1}^U \neq \emptyset$,
then
$$
N_{i+1}^U -U\subseteq \cup_{V\subset U} N_{i}^V  =N^{\subseteq U}_{i}.
$$
By the Cauchy-Davenport theorem, $|N_{i+1}^U -U|\ge |N_{i+1}^U |+|U|-1$.
\end{proof}

We will use the following result originally obtained in \cite{SZ}.
We provide here a shortened proof. We use the notation $\overline{X}$
as shorthand for $\zp\setminus X$.

\begin{theorem}{\rm (\cite{SZ})}
Let  $B$  be a  subset of
 $\zp$ containing $0$ and let $A$ be a $2$-atom of $B$ containing $0$.
Set $m=\kappa_2(B)-|B|=|B+A|-|B|-|A|$.
Assume  that  $|B|<p-(m+4)(m+3)/2$:
then $|A|=2.$ In particular $B$ is a union of at most $m+2$ arithmetic
progressions with the same difference.
 \label{SZ}
 \end{theorem}

 \begin{proof}
   Suppose that $|A|>2$. Set $A^*=A\setminus \{ 0\}.$
Let us first show that
\begin{equation}
N_{i+1}-A^*\subset N_{i} , \;i\ge 1, \label{eq:a*}
\end{equation}
where $N_i=N_i(B,A)$. Take $x \in N_2$. Since $x-A$ meets $N_1$, then
$x-A$ is not a
 subset of $\overline{B+A}$. By Lemma \ref{dual}, $x-A$ and
 $\overline{B+A}$ are respectively a $2$-atom and a $2$-fragment of
 $-B$. Theorem \ref{intersection} implies
 $(x-A)\cap \overline{B+A}=\{x\}$ or, equivalently, $x-A^*\subset
 N_1$. We have just proved   \eqref{eq:a*} for $i=1$.
 The assertion follows for $i>1$ by Lemma \ref{lem:NU} since $N_2^{U}=\emptyset$ for
 each proper subset $U\subset A^*$.

Let $t$ be the largest integer such that $N_t\neq \emptyset$. Since $A$
generates $\zp$, we have
$$
\zp=B\cup N_1\cup N_2 \cup \cdots \cup N_t, \; t\ge 2,
$$

By \eqref{eq:a*} and Lemma \ref{lem:NU} we have
\begin{equation}
|N_{i+1}|\leq |N_{i}|-(|A|-2),\; 1\le i <t, \label{eq:ni}
\end{equation}
and  $|N_1|=|B+A|-|B|=|A|+m$.

Suppose first that $t\ge 3$. Using
\eqref{eq:ni} and $|A|\ge 3$,
\begin{eqnarray*}
p&\le&|B|+(m+|A|)+(m+2)+(m+4-|A|)+\sum_{i=0}^{m-1} (m-i) \\
&\le &|B|+\sum_{i=1}^{m+3}i =|B|+(m+4)(m+3)/2,
\end{eqnarray*}
contradicting the assumption on $|B|$.

Suppose now that $t=2$. In this case $\overline{B+A}=N_2$.
Since
$\overline{B+A}$ is a $2$--fragment of $-B$, we have $|A|\le
|-(\overline{B+A})|=|N_2|$: furthermore, $|N_2|\le m+2$ by
\eqref{eq:ni}.
We obtain therefore $p\le |B|+(m+|A|)+(m+2)<|B|+(m+4)(m+3)/2$,
 again a contradiction.
 \end{proof}

Theorem \ref{SZ} provides us with a useful bound of the size of $k$-atoms.

\begin{proposition}\label{prop:size}
Let  $B$  be a  subset of
 $\zp$ containing $0$ and let $A$ be a $k$-atom of $B$ with $2\leq k\leq |B|$.
 Put $m=\kappa _k(B)-|B|$.
 Assume moreover $p+k> m^2+6m+12.$
 Then $|A|\leq m+k+1.$
\end{proposition}

\begin{proof} Without loss of generality, we may assume $0\in A$.
By definition of $m$ we have $|B+A|=|B|+|A|+m$. Let
$C=\overline{B+A}$ so that we have $p=|A|+|B|+|C|+m$. By
definition of $k$-atoms we have $|C|\geq k$ and $|A|\geq k$.
Observe that $(C-B)\cap A=\emptyset$ so that $|B-C|=|C-B|\leq
|B|+|C|+m$ and therefore $|A|\leq |C|$ by the minimality condition
in the definition of $k$-atoms. Hence, $2|A|\leq |A|+|C|
=p-|B|-m$. The condition $k\leq |B|$ now implies
$$
|A|\leq (p-k-m)/2<
(2p-m^2-6m-12-m)/2=p-(m+4)(m+3)/2.
$$
Now $|B+A|=|B|+|A|+m$ implies that $\kappa_k(A)-|A|\leq m$ and hence
that $\kappa_2(A)-|A|\leq m$. Therefore Theorem \ref{SZ} implies that
$A$ is a union of not more than $m+2$ arithmetic progressions with the
same difference $u$. This gives
  $$|A|-m-2\le |A\cap (A+u)|.$$
 Furthermore, since $A+u$ is also
a $k$-atom of $B$, Theorem \ref{intersection} implies
  $$|A\cap (A+u)|\le k-1,$$
hence the result.
\end{proof}

\section{Compression transfer}\label{sec:comp}

Let us start with a lemma which is the $\z$ counterpart of our main result.

\begin{lemma}
\label{lem:zcase}
Let $A$ and $B$ be subsets of $\z$ such that $0\in A\cap B$ and
$$
|A+B|\leq |A|+|B|+1,\;\; |B|\ge |A|\ge 4, |B|\ge 5.
$$
Then $A$ and $B$ have the same greatest common divisor $r=\gcd(A)=\gcd(B)$
and  $\ell_r(A)\leq |A|+2$ and $\ell_r(B)\leq |B|+2$.
\end{lemma}

\begin{proof}
Let $r=\gcd(A)$. Put $B_1=\{x\in B,\; x=0\;\bmod r\}$
and $B_2=B\setminus B_1$.
We have $B_1\neq \emptyset$ since $0\in B_1$ and if
$\gcd(B) < r$ then $B_2\neq \emptyset$ also. We then must
have
$|A+B|\geq |A|+|B_1|-1+|A|+|B_2|-1\geq |A|+|B|+2$, a contradiction.
This proves $\gcd(B)\geq\gcd(A)$. Proceed likewise to obtain
$\gcd(A)\geq\gcd(B)$.

We may assume
$A,B\subset \n$ and $\gcd(A)= \gcd(B)=1$ without loss of generality.
Put $a=\max (A)$ and $b=\max (B)$.

Assume first $a\geq b$.
Apply Theorem \ref{thm:ls} with $A$ and $B$ interchanged
to get
  $$|A+B|\geq \min(a+|B|,2|B|-1 +|A|-2).$$
Hence $|A|+|B|+1\geq \min(a+|B|,|A|+|B|+(|B|-3))$,
and since $|B|\geq 5$ we must have $|A|+|B|+1\geq a+|B|$,
which means $\ell_1(A)\leq |A|+2$. Since $|A|\leq |B|$ and
$\max(B)\leq \max(A)$ we
clearly must also have $\ell_1(B)\leq |B|+2$.

Suppose now $a<b$. Apply Theorem \ref{thm:ls} to get
  $$|A+B|\geq \min(b+|A|,2|A| +|B|-2).$$
 Since $|A|\geq 4$ we get $|A|+|B|+1 \geq b+|A|$
which means $\ell_1(B)\leq |B|+2$. Notice furthermore that if
$|B|\leq |A|+1$ then $a<b$ implies $\ell_1(A)\leq |A|+2$.

Summarizing, we have proved that under the hypothesis of the lemma,
$\ell_1(B)\leq |B|+2$ always holds and $\ell_1(A)\leq |A|+2$ holds
under the additional condition
\begin{equation}
  \label{eq:additional}
  |B|\leq |A|+1.
\end{equation}
We proceed to prove by induction on $|B|$ that under the hypothesis
of the lemma $\ell_1(A)\leq |A|+2$
always holds. If $|B|=5$ then \eqref{eq:additional} must hold and we
are done.
Now suppose by induction that the result holds for $|B|=\beta \geq 5$, and
consider the case $|B|=\beta +1$.
Let $B'=B\setminus\{b\}$. Since the result holds with \eqref{eq:additional}
assume $|A|\leq |B'|$. Since $|B|=\beta +1\geq 6$
we have $|B'|\geq 5$.
Since
$a+b\not\in A+B'$ we also have $|A+B'|\leq |A|+|B'|+1$.
The result follows from the induction hypothesis applied to $A$ and $B'$.
\end{proof}

The aim of this section is to show that, under the hypothesis of Theorem
\ref{thm:main}, if one of the two sets is in a short arithmetic progression
with difference $r$, then so is the other one.

\begin{theorem} \label{thm:collage}
Let $X$, $Y$ be subsets of $\zp$ such that
$$
|X+Y|=|X|+|Y|+1 \leq p-5,
$$
and with $|Y|\ge 4$, $|X|\ge 5$.
If $p> 32$
then $\ell_1(Y)\leq |Y|+2$ implies $\ell_1(X) \leq |X|+2$.
\end{theorem}

The proof of Theorem \ref{thm:collage} will be broken down into
several lemmas.

First, we need some notation.
By a {\it connected component} of a set $Z\subset \zp$ we mean a maximal arithmetic progression of difference $1$ contained in $Z$. Let $C_1,\ldots ,C_j$ be
the connected components of the complement $\ox$ of $X$ with
$|C_j|=\max_{1\le i\le j}|C_i|$.
Thus, $\ell_1(X)=p-|C_j|$. We have
\begin{equation}\label{eq:ci}
|X|+|Y|+1= |X+Y|=|X|+\sum_{i=1}^j|(X+Y)\cap C_i|.
\end{equation}
For $i=1,\ldots ,j$ we shall use the notation
$C_i=\{ c_i, \ldots ,c_i+ |C_i|-1\}$ and   $X_i$   denotes the connected component of $X$ containing $c_i-1$. We  assume that $j>1$ since otherwise there is nothing to prove.
\begin{lemma}\label{bigcfr} Theorem \ref{thm:collage} holds if $|C_j|\ge \ell_1(Y)-1$.
\end{lemma}
\begin{proof}
$|C_j|\ge \ell_1(Y)-1$ means that by translation we may choose
$Y\subset [0,\ell_1(Y)-1]$ and $X\subset [0,\ell_1(X)-1]$ with
$\ell_1(X)+\ell_1(Y)-1\leq p$. In other words $X+Y$ can really be
considered as a sum in $\z$ and the result follows from Lemma
\ref{lem:zcase}.
\end{proof}
\begin{lemma}\label{lem:tap+1} Theorem \ref{thm:collage} holds if
 $|C_j|< \ell_1(Y)-1$ and $\ell_1(Y)\le |Y|+1$.
\end{lemma}
\begin{proof}
Suppose that $\ell_1(Y)=|Y|$. Then $|(X+Y)\cap C_i|=|C_i|$ which
implies $X+Y=\zp$ against the assumptions.

Suppose now that $\ell_1(Y)=|Y|+1$.
Let $Y_1, Y_2$ be the connected components of $Y$ with $y_i=|Y_i|$. By using multiplication by $-1$ and translating if necessary, we may assume that  $Y_1=\{ 0,\ldots ,y_1-1\}$, $Y_2=\{ y_1+1,\ldots ,y_1+y_2\}$ and $y_1\ge y_2$.

We have
\begin{equation}\label{eq:shortc}
 |(X+Y)\cap C_i|=\left\{\begin{array}{ll}|C_i|, & \mbox{ or }\\ |C_i|-1\ge |Y_1|-1,& \mbox{  and  } c_i-2\in \ox,\end{array}\right.
\end{equation}
which implies $j\ge |\oxy|$ and, by using  (\ref{eq:ci}),
\begin{eqnarray*}
|Y|+1&=&\sum_{i=1}^j|(X+Y)\cap C_i|\ge |\oxy|(|Y_1|-1)\ge 5(|Y_1|-1)
\\ &\ge& |Y|+3|Y_1|-5.
\end{eqnarray*}
Hence $|Y_1|=2$ and all equalities hold. In particular, $X$ consists
of $|\oxy|=5$ connected components each of cardinality one and
$|C_i|=2$ for each $i$, which implies that $p$ is divisible by $3$, a
contradiction.
\end{proof}

In the remaining part of this section we  assume that
$\ell_1(Y)=|Y|+2$. 
Let $Y=Y_1\cup Y_2\cup Y_3$ be the decomposition of $Y$ into connected
components, with $Y_2$ possibly empty,    $y_i=|Y_i|$ and $y_1=\max_{i} y_i$.
By using multiplication by $-1$ and/or translating $Y$ if necessary,
we shall assume that $Y_1\cup Y_2\subset \{ 0,1,\ldots ,y_1+y_2 \}$.
\begin{lemma}\label{lem:t>12}
  Theorem \ref{thm:collage} holds if $|C_j|<\ell_1(Y)-1$ and
either $|Y|>9$ or $|\oxy|\geq 12$.
\end{lemma}
\begin{proof} Since $|C_i|<\ell_1(Y)-1$ we have $|(X_i+Y)\cap C_i|\ge |C_i|-2$ and, if equality holds, then
$|C_i|\ge |((c_i-1)+Y)\cap C_i|\ge |Y_1|+|Y_2|-1$. Moreover, $c_i-2\in \ox$. By using similar remarks when $|(X_i+Y)\cap C_i|= |C_i|-1$ we have
\begin{equation}\label{eq:ci2}
|(X+Y)\cap C_i|=\left\{\begin{array}{ll} |C_i|&\mbox{  or  }\\
|C_i|-1\ge \lceil |Y|/3\rceil-1 &\mbox{  and  } |X_i|\le 2, \mbox{ or}\\
|C_i|-2\ge   \lceil |Y|/2\rceil-1&\mbox{  and  } |X_i|=1.
\end{array}\right.
\end{equation}

Suppose first that $|Y|\ge 10$. For $r=0,1,2$ let $J_r\subset J=\{ 1,\ldots ,j\}$ such that
$|(X+Y)\cap C_i|=|C_i|-r$.  We have
$|J_1|+2|J_2|=|\ost|\ge 5$. Therefore, using (\ref{eq:ci}) and (\ref{eq:ci2}),
$$
|Y|+1=\sum_{i=1}^j|(X+Y)\cap C_i|\ge |J_0|+|J_1|(\lceil |Y|/3\rceil-1)+|J_2|(\lceil |Y|/2\rceil-1),
$$
which can hold only if $|J_0|=0$, $|J_1|=1$ and $|J_2|=2$ (and $|Y|=10$ or $12$.) In this case (\ref{eq:ci2}) implies that $X$ consists
of two connected components of cardinality one and one connected component of cardinality at most two, against the assumption that $|X|\ge 5$.

Suppose now that  $|Y|\leq 9$ and $|\oxy|\geq 12$.
Then (\ref{eq:ci2}) implies that the connected components of $\oxy$ have
cardinality $1$ or $2$.
If at most one has cardinality $2$, then
$|J|>10$. On the other hand, since   $|Y|\geq 4$ we have $|Y_1|\ge 2$ and $|(X+Y)\cap C_i|\ge 1$ for each $i\in J$. Hence,
$|X+Y|>|X|+10$ which contradicts $|X+Y|=|X|+|Y|+1$.
If at least two connected components of $\oxy$ have cardinality $2$,
then (\ref{eq:ci2}) implies that the corresponding two $C_i$'s contain at
least $|Y|-2$ elements of $(X+Y)\setminus X$. The remaining $C_i$'s
contribute at least one element to $(X+Y)\setminus X$ and there are at least
(12-4)/2=4 of them, meaning again that $|X+Y|> |X|+|Y|+1$.
\end{proof}

\medskip

\noindent
{\bf Proof of Theorem \ref{thm:collage}.}
By Lemmas \ref{bigcfr}, \ref{lem:tap+1} and \ref{lem:t>12}, the only case left to be
examined is when $\ell_1(Y)=|Y|+2$, $|Y|\leq 9$ and $|\oxy|\leq 11.$ Let us set $X' = \oxy$.
Note that we have $X'-Y \subset \ox$ so that $|X'-Y| \le |X'| + |-Y| +1$. Now
the condition $p>32$ implies  $|\overline{X'-Y}|=p-|X'|-|Y|-1>11$. Therefore Lemma \ref{lem:t>12} applies to $X'$ and $-Y$
and we obtain that  $\ell_1(X')\leq |X'|+2$. Hence
$X'-Y=\ox$ also satisfies $\ell_1(\ox)\leq |\ox|+2$ (easily checked).
But this means that
$X=\overline{\ox}$   is the union  of at most two single elements
and of a progression $Z$ with $|Z|\ge 3$, since $|X|\geq 5$.
By (\ref{eq:ci2}) $|Z|\geq 3$ implies $|J_0|\geq 1$ and
since $|J|\leq 3$ we get $|\oxy|\le 4$, a contradiction.
This concludes the proof of
Theorem \ref{thm:collage}.

\section{On the size of atoms}\label{sec:sizebis}

The bound on the size of  $k$-atoms of a set $X$ given by Proposition
 \ref{prop:size} can be improved when $\kappa_k (X)\le |X|+1$. 
We prove the following.

\begin{theorem}\label{thm:size2}
Let $X\subset \zp$ such that $0\in X$, $|X|\ge 4$ and $\kappa_5(X)\le |X|+1$.
Suppose $p> |X|+42$. Then 
\begin{itemize}
\item[(i)]  the $4$-atoms of $X$ have cardinality $4$,
\item[(ii)] if $|X|=4$ then the $5$-atoms of $X$ have cardinality
$5$.
\end{itemize}
\end{theorem}

We shall break up the proof of Theorem \ref{thm:size2} into several lemmas.
First we introduce some notation and terminology that will be
convenient to us in this section.

For $z\in N_i(X,Y)$, $i\geq 1$, define its {\em outdegree}
$d_+(z) = |(z+Y)\cap N_{i+1}|$ and its {\em indegree}
$d_-(z) = |(z-Y)\cap N_{i-1}|$. Note that
by counting in two ways the number of couples $(z,z')$ such that
$z\in N_i$, $z'\in N_{i+1}$ and $z'-z\in Y$, we have~:

\begin{equation}\label{eq:balance}
\sum_{z\in N_i}d_+(z) = \sum_{z\in N_{i+1}}d_-(z).
\end{equation}

We shall call the quantity in \eqref{eq:balance} indifferently the
total outdegree of $N_i$ or the total indegree of $N_{i+1}$.

\begin{lemma}\label{lem:n1} 
Suppose $A$ is a $k$-atom of some set $X\subset\zp$.
If there is $z\in X+A$
uniquely expressable as $z=x+a$, $x\in X$ and $a\in A$ then
$|A|=k$.
\end{lemma}
\begin{proof} If $|A|>k$ then $A'=A\setminus \{ a\}$ satisfies
$|X+A'| -|A'| \leq |X+A|-|A|$ which contradicts $A$ being a
$k$-atom.
\end{proof}

\begin{lemma}\label{lem:rect1} Theorem \ref{thm:size2} holds if
$\ell_r(X)\le |X|+2$ for some $r$.
\end{lemma}
\begin{proof} Suppose on the contrary that $A$ is a $k$-atom of $X$ of size
$|A|>k\ge 4$. By Theorem \ref{thm:collage} we have  $\ell_r(A)\le
|A|+2$ so that the sum $X+A$ can be considered as a sum in $\z$.
There is therefore  $z\in X+A$ uniquely expressable as $z=x+a$,
$x\in X$, $a\in A$, contradicting Lemma \ref{lem:n1}.
\end{proof}

\begin{lemma} \label{lem:y=6}
Let $X$ and $p$ be as in Theorem \ref{thm:size2}. Let $A$ be a
$4$-atom of $X$. Then $|A|\le 5$.
\end{lemma}

\begin{proof} By Proposition \ref{prop:size} we have $|A|\le 6$.
Suppose   that $|A|=6$.  We may assume that $0\in A\cap X$. Let
$A^*=A\setminus \{ 0\}$. Note that $A^*$ is not a $d$-progression,
since otherwise $|A\cap (A+d)|\ge 4$ contradicting 
Theorem~\ref{intersection}. By Lemma \ref{lem:rect1} we may assume $\ell_r
(X)\ge |X|+3$ for each $r\in \zp^*$.
\begin{enumerate}
\item
We have $|N_2^{ab}|=0$ for every $a,b\in A^*$.

Otherwise $|N_1^{ab}|\geq 2$ and $A'=A\setminus\{a,b\}$ satisfies
$|X+A'| -|A'| \leq |X+A|-|A|$ which contradicts $A$ being a
$4$-atom.
\item
We have $|N_2^{abc}|\leq 1$ for every $a,b,c\in A^*$.

Otherwise $|N_1^{\leq abc}|\geq 4$ and $A'=A\setminus\{a,b,c\}$
satisfies $|(X+A')\setminus X| \leq 3$ but then Theorem
\ref{thm:hamrod} implies that $\ell_r(X)\leq |X|+1$ for some $r$.

Lemma \ref{lem:NU} implies therefore
\item
$|N_3^{abc}| = 0$ for every $a,b,c\in A^*$.
\item
We have $|N_2^{abcd}|\leq 2$ for every $a,b,c,d\in A^*$.

Otherwise $|N_1^{\leq abcd}|\geq 6$ and
$\{0,e\}=A\setminus\{a,b,c,d\}$ is such that $|X+\{0,e\}| \leq
|X|+1$ which implies that $\ell_e(X)=|X|$.
\item
We have $|N_3^{abcd}|\leq 3$ for every $a,b,c,d\in A^*$.

Otherwise $|N_2^{\leq abcd}|\geq 7$ by Lemma \ref{lem:NU}, but
this contradicts $|N_2^{abcd}|\leq 2$ and $|N_2^U|\leq 1$ for
$|U|=3$. Hence we get
\item
$|N_4^{abcd}|=0$ for every $a,b,c,d\in A^*$, i.e. $N_i=N_i^{A^*}$
for $i\geq 4$.
\end{enumerate}
We now bound from above the $|N_i|$. Since $|N_1|\leq 7$, the
total outdegree of $N_1$ is at most $7\times 5 =35$, which
implies, since the indegree of any element of $N_2$ is at least
$3$, that
  $$|N_2|\leq 11.$$

This means that the total outdegree of $N_2$ is at most $11\times
5 = 55$, which in turn implies, since the indegree of any element
of $N_3$ is at least $4$, that
  $$|N_3|\leq 13.$$
Since $N_i=N_i^{A^*}$ for $i\ge 4$ and $A^*$ is not an arithmetic
progression, we have, by Vosper's Theorem,
$|N_i|+|A^*|\le |N_i-A^*|\le |N_{i-1}|$ which
implies $|N_4|\leq 8$ and $|N_5|\leq 3$. This implies
$|\overline{X}|\leq 42$.
\end{proof}

\begin{lemma} \label{lem:y=5}
Let $X$ and $p$ be as in Theorem \ref{thm:size2}. Let $A$ be a
$4$-atom of $X$. Then $|A|=4$.
\end{lemma}

\begin{proof}
Suppose the contrary. We may assume $0\in A$.
By Lemma \ref{lem:y=6} we may also assume $|A|=5$. 
First note that $A$ cannot be an arithmetic progression of some
difference $d$, since $|A\cap (A+d)|=4$ would contradict 
Theorem~\ref{intersection}. 
 Note that, if $A^*$ is
an arithmetic progression, then $A$ is equivalent to
$B=\{0,1,2,3,u\}$. If $B^*$ were an arithmetic progression
of difference $d$ we would have $|B^*\cap (B^*+d)|=3$ and
$|\{1,2,3\}\cap\{1+d,2+d,3+d\}|\geq 1$ implying that
$B^*$ is an arithmetic progression of difference $1$ or $2$,
none of which are possible. We may therefore assume,
without loss of generality, that $A$ is a $4$-atom of $X$
containing $0$ such that  $A^*$ is not an arithmetic progression.

By Lemma \ref{lem:rect1} we may also assume that
$\ell_r(X)\ge |X|+3$ for each $r\in \zp^*$.
We proceed very much
along the same lines as in the previous Lemma.

\begin{enumerate}
\item We have $|N_1^{ab}|\leq 2$ for each $a,b\in A^*$.

Otherwise $|X+(A\setminus \{ a,b\})|\leq |X|+3$ and Theorem
\ref{thm:hamrod} implies that $\ell_r(X)\leq |X|+1$ for some $r$.

By Lemma \ref{lem:NU} we have

\item $|N_2^{ab}|\leq 1$ and $|N_i^{ab}|=0$, $i\ge 3$, for each $a,b\in A^*$.

\item $|N_2^{abc}|\leq 2$ for each $a,b,c\in A^*$.

Otherwise, with $V=\{ a,b,c\}$, we have  $|N_2^{V}-V|\geq 5$ and
$|X+ (A\setminus V)| \leq |X|+1$ implying that $X$ is an
arithmetic progression.

Lemma \ref{lem:NU} then implies $|N_3^{V}-V|\le  |N_2^V| +
|N_2^{ab}|+|N_2^{ac}|+|N_2^{bc}|\le 5$, so that we get
\item $|N_3^{abc}|\le 3$ and $|N_4^{abc}|\leq 1$ for each $a,b,c\in A^*$.

In particular,

\item $N_i=N_i^{A^*}$ for $i\geq 5$.
\end{enumerate}

We now bound from above the number of elements in $N_i$, $i\geq
2$. Since $|N_1|\leq 6$ and $|N_2^{ab}|= 1$ implies
$|N_1^{ab}|=2$,  there are at most $3$ elements of $N_2$ of
indegree $2$. Now since the total outdegree of $N_1$ is at most
$|N_1|\cdot|A^*|=6\times 4= 24$, we get by (\ref{eq:balance}) that
$|N_2|\leq 9$. Actually we must have
  $$|N_2|\leq 8,$$
because $|N_2|=9$ can occur only if  $d_+(z)=4$ for every $z\in
N_1$ which implies $|N_2|= |N_1 + A^*|=|N_1|+|A^*|-1$ and, by
Vosper's theorem,  $A^*$ is an arithmetic progression against our
assumption.

Now the total outdegree of $N_2$ is at most $8\times 4 =32$, and
every element of $N_3$ has indegree at least $3$, so that
(\ref{eq:balance}) implies
  $$|N_3|\leq 10.$$
Finally, $|N_4^{A^*}|+|A^*|\le |N_4^{A^*}-A^*|\le |N_3|$ implies
$|N_4^{A^*}|\le 6$. Since the sets $N_3^V$ for $|V|=3$ and
$N_4^V\neq \emptyset$ are disjoint and contain at least $3$
elements, we have $|N_4\setminus N_4^{A^*}|\leq 3$ and
 $$|N_4|\leq 9.$$
For every $i\ge 5$ we have $|N_i|=|N_i^{A^*}|\le |N_{i-1}|-|A^*|$
so that $|N_5|\le 5$ and $|N_6|\le 1$. Adding up the $N_i$'s we
get $|\overline{X}|\le 39$.
\end{proof}

Lemma \ref{lem:y=5} proves point $(i)$ of Theorem
\ref{thm:size2}. To prove point $(ii)$
we use the following Lemma.

\begin{lemma}\label{lem:y=4}
Let $X,B\subset \zp$,  such that $|X|= 4$, $|B|\geq 4$, 
$0\in X\cap B$, and
  $$|X+B|\leq |X|+|B|+1.$$
If $p>|B|+20$, then there is an element $z\in X+B$ which can be
uniquely written as $z=x+b$ with $x\in X$ and $b\in B$.
\end{lemma}

\begin{proof}
If $B$ is an arithmetic progression of difference $d$ then 
$|X+B|\leq |X|+|B|+1$ clearly implies $\ell_d(X)\leq |X|+2$,
so that the sum $X+B$ can be considered as a sum in $\z$, in which
case the conclusion of the lemma holds. Suppose therefore that
$B$ is not an arithmetic progression. 
Let $X^*=X\setminus \{ 0\}$. We now write $N_i=N_i(B,X)$.
Suppose that there does not exist $z$ uniquely expressable as
$z=x+b$, $x\in X$, $b\in B$. This implies
$N_1^{x}=\emptyset$ for every $x\in X^*$.

We have $|N_1^{xy}|\le 3$ for each $x,y\in X^*$.

Otherwise $|B+(X\setminus \{x,y\})|\le |B+X|-|N_1^{xy}|\le |B|+1$
implying that $B$ is an arithmetic progression.

By successively applying Lemma \ref{eq1} and writing $(m)^+=\max\{
0,m\}$, we have
\begin{eqnarray*}
|N_2|&=&\sum_{V\subset X^*}|N_2^V| = \sum_{V\subset X^*,
|V|=2}|N_2^V|
         + |N_2^{X^*}|\\
 &\le &\sum_{V\subset X^*,
|V|=2}(|N_1^V|-1)^+ +(|N_1|-2)\le 6,
\end{eqnarray*}
where the last inequality uses the fact that $|N_1^V|=3$ and
$|V|=2$ occur together at most once. Similarly,
\begin{eqnarray*}
|N_3|&\le& \sum_{V\subset X^*,
|V|=2}(|N_2^V|-1)^++(|N_2|-2)\le 5, \\
|N_4|&=&|N_4^{X^*}|\le |N_3|-2\le 3,\;\;\mbox {and }\\
|N_5|&\le &|N_4|-2\le 1.
\end{eqnarray*}

Therefore, $p\le |B|+\sum_{i=1}^5|N_i|\le |B|+20$.
\end{proof}

To prove point $(ii)$ of Theorem \ref{thm:size2} consider
a $5$-atom $B$ of $X$. By Proposition~\ref{prop:size} we have $|B|\leq 7$,
and since $p\geq 29$ Lemma \ref{lem:y=4} implies that there is $z$
uniquely expressable as $z=x+b$, $x\in X$, $b\in B$. But this
contradicts Lemma \ref{lem:n1}.

\section{The case of small sets}\label{sec:small}

We next prove Theorem \ref{thm:main} when the two sets attain
their minimum possible values, $|A|=4$ and $|B|=5$. We first need
the following two lemmas. As in section \ref{sec:comp} we call
a $d$--component of a set
$Z\subset \zp$ a maximal arithmetic progression of
difference $d$ contained in $Z$. We denote by   $c_d(Z)$
 the number of $d$-components of $Z$.

 \begin{lemma}\label{lem:2comp}
  Let $p>23$ and let $A,B\subset \zp$ with $|A|=4$, $|B|=5$ and
$|A+B|=10$. Then $c_d(A)\le 2$ for some $d\in \zp$.
\end{lemma}

\begin{proof}
  Without loss of generality we may suppose $0\in A\cap B$.
  By Lemma~\ref{lem:y=4}   there exists $C\subset B$,
  $|C|=4$, such that $|A+C|<|A+B|$. By inclusion-exclusion,
$$
9\ge  |A+C|=|\cup_{c\in C}(A+c)|\ge 16-\sum_{c, c'\in C, c\neq
c'}|(A+c)\cap (A+c')|,
$$
which gives $\sum_{c, c'\in C, c\neq c'}|(A+c)\cap (A+c')|\ge
7>\binom{|C|}{2}$. Hence there are two distinct elements  $c,
c'\in C$ such that $|(A+c)\cap (A+c')|\ge 2$   or, equivalently,
that $|A+\{0,c-c'\}|\le |A|+2$. It follows that $C$ has at most
two $d$-components for $d=c-c'$.\end{proof}

\begin{lemma}\label{lem:rect}
Let $p$ be any odd prime and let $Z\subset\zp$ with $0\in Z$ and
$|Z|< (p+9)/4$.  If $c_d(Z)\le 2$  for some $d\in (\zp)^*$
then some affine image of $Z$ is a subset of $ \{0, 1, \cdots ,
(p-1)/2\}$.
\end{lemma}

\begin{proof} Suppose  that no  affine image of $Z$ is  a subset
of $ \{0, 1, \cdots , (p-1)/2\}$. Then we must have
$c_d(Z)=2$.

Let $A$ and $B$ be the components of $Z$, where $a=|A|\geq |B|=b$.
Let $V$ and $W$ be the components of $\overline{Z}$, where
$v=|V|\leq |W|$. Without loss of generality   we may assume that
$A$ is represented by the integers $A_0=\{0, 1, \cdots , a-1\}$
and $B$ is represented by $B_0=\{v+a,  \cdots , v+a+b-1\}$.

Since  $Z\not\subset \{0, 1, \cdots , (p-1)/2\}$, we have
$v+a+b-1\geq (p+1)/2$. Moreover $2v+a+b\leq p$. It follows that
$(p+1)/2-a-b+1\leq v \leq  (p-a-b)/2$.
Now $ B_0\subset [(p+1)/2-b+1,(p+a+b)/2-1 ]$. It follows that
 $ 2*B_0\subset [p+3-2b, p+a+b-2 ]$ so that $2*B$ is represented by a
subset of $\{-2b+3,-2b+2,\ldots , a+b-2\}$. Since $2*A$ is represented
by a subset of $\{0,2,\ldots , 2a\}$, we get that
$2*Z$ is represented by a subset of
$\{\inf\{-2b+3,0\}, \ldots , \sup\{a+b-2,2a-2\}\}$. Since $a\geq b$ we
get that $2*Z$ is represented by a subset of
$\{-2b+2,\ldots ,2a-2\}$. Now we must have $(p+1)/2\leq
2a+2b-4=2|Z|-4$, and hence $p\leq 4|Z|-9$, a contradiction.
\end{proof}

\begin{lemma}\label{lem:45}
Let   $p> 23$ and $A, B\subset \zp$ with   $0\in A\cap B$, $|A|=4$
and $|B|=5$ and $|A+B|=10$.

Then there is $r\in \zp$ such that $\ell_r (A) \leq |A|+2$ and
$\ell_r (B) \leq |B|+2$.
\end{lemma}

\begin{proof}
Asume first that $c_d(A+B)\leq 2$ for some $d\in\zp$. Then by
Lemma \ref{lem:rect}, $Z=A+B$ is such that $\ell_a(Z)\leq (p-1)/2$
for some $a$.
This implies that if $z_1,z_2,z_1',z_2'$ are four integers
representing elements of $Z$, then $z_1+z_2=z_1'+z_2' \bmod p$
implies $z_1+z_2=z_1'+z_2'$. Since $0\in A\cap B$ we have
$A\cup B\subset Z$ therefore $A+B$ can be considered as a sum
in $\z$ and Lemma \ref{lem:zcase} implies the result.

Assume now that $A+B$ has at least three $x$-components for every $x\in\zp$.

By Lemma \ref{lem:2comp} we have $c_d(A)\le 2$ for some $d\in\zp$.
Assume $c_d(A)=2$ otherwise there is nothing to prove.
Let us show that we must have both
\begin{equation}\label{claim:u}
  c_d(A)=2,\hspace{1cm}c_d(B)\le 2.
\end{equation}
Suppose the contrary. Set $A=A_1\cup A_2$ be the decomposition of
$A$ into $d$-components where $|A_1|\geq |A_2|$.

Let us first show that $|A_1|=2$. Suppose the contrary, i.e. $|A_1|=3$. Since
$B$ has at least three $d$-components, we have
$|A_1+B|\geq |B|+3+t$, where $t$ is the number of $d$-components
of $A_1+B$. Since $A+B$ also has at least three $d$-components, we have
$|A+B|\geq |B|+3+t+(3-t)=11$, a contradiction.

Now we can write $A=\{0,d\}\cup \{x,x+d\}=\{0,d\}+\{0,x\}$.
Since $B$ has at least three $d$-components we have
$|B+\{0,d\}|\geq |B|+3$. Now observe that $B+\{0,d\}$
has at least three $x$-components since otherwise
$A+B=(B+\{0,d\})+\{ 0,x\}$ would have less than three
$x$-components. It follows that $|A+B|=|B+\{0,d\}+\{0,x\}|\geq
|B|+6$, a contradiction. This proves~\eqref{claim:u}.

We may now assume $c_d(B)>1$, i.e. $c_d(B)=2$, otherwise there is nothing
to prove. Let $B=B_1\cup B_2$ be the $d$-components of $B$.
 We have
$c_d(A+B)=3$, since otherwise $A+B$ has 4 components that are
necessarily $A_1+B_1, A_1+B_2, A_2+B_1, A_2+B_2$. Since
$|A_i+B_j|=|A_i|+|B_j|-1$ we would get $|A+B|=2|A|+2|B|-4=14$.

Observe that since $A_1+B_1$ and $A_1+B_2$ are disjoint (because
$|A_1|=2$), they must belong to distinct
$d$-components, otherwise all the sets $A_i+B_j$ are disjoint and
we again get $|A+B|=14$. For the same reason $A_2+B_1$ and $A_2+B_2$
must belong to distinct $d$-components.

The sets $A_1+B_1$ and $A_2+B_1$ also belong to distinct
$d$-components, otherwise this common component $C$ satisfies $|C|\geq
|A|+1+|B_1|-1= |A|+|B_1|$. Then the other components are $A_1+B_2$
and $A_2+B_2$. It follows that  $|A+B|\geq
2|A|+|B_1|+2|B_2|-2>10$. For the same reason $A_1+B_2$ and $A_2+B_2$
belong to distinct $d$-components.

The only possibility left is that $A_1+B_1\cup A_2+B_2$ merge into one
component $C$. Since one of the two $d$-components of $B$ has at least
three elements we have $|C|\geq 4$. The remaining components of $A+B$
have $|A_1|+|B_2|-1$ and $|A_2|+|B_1|-1$ elements which gives
$|A+B|\geq 4 + |A|+|B|-2=11$, again a contradiction. This completes
the proof.
\end{proof}

\section{Proof of Theorem \ref{thm:main}}\label{sec:main}

We are now ready for the proof of Theorem \ref{thm:main}

Suppose that $|A|=4$. Let $U$ be a $5$-atom of $A$. By Theorem
\ref{thm:size2} we have $|U|=5$. By Lemma \ref{lem:45} we have
$\ell_r(A)\le |A|+2$   for some $r\in (\zp)^*$ and by Theorem
\ref{thm:collage} we then have $\ell_r(B)\le |B|+2$.

Suppose now that $|A|\ge 5$. Let $U$ be a $4$-atom of $A$.
By proposition \ref{prop:size} we have $|U|\leq 6$.

\begin{enumerate}
\item
If $|U|=4$. Then let $V$ be a $5$-atom of $U$. By Theorem
\ref{thm:size2} we have $|V|=5$, and by Lemma \ref{lem:45} we have
$\ell_r(U)\le |U|+2$ for some $r\in (\zp)^*$. By Theorem
\ref{thm:collage} we then have $\ell_r(A)\le |A|+2$, and again by
Theorem \ref{thm:collage} we finally have $\ell_r(B)\le |B|+2$.
\item
If $5\leq |U|\leq 6$. Then let $V$ be a $4$-atom of $U$. We have
$p\geq 53$ which implies $p> |U|+42$ and therefore Theorem
\ref{thm:size2} implies $|V|=4$. Let $W$ be a $5$-atom of $V$.
Apply Theorem \ref{thm:size2} again to obtain $|W|=5$. Then Lemma
\ref{lem:45} implies $\ell_r(V)\le |V|+2$ for some $r\in (\zp)^*$.
Theorem \ref{thm:collage} applied once implies that $\ell_r(U)\le
|U|+2$ then applied once more implies that $\ell_r(A)\le |A|+2$
and applied a third time implies finally that $\ell_r(B)\le
|B|+2$. This concludes the proof.
\end{enumerate}

\end{document}